# Flying Higher Than A Box-Kite:

# Kite-Chain Middens, Sand Mandalas, and Zero-Divisor Patterns in the $2^n$-ions Beyond the Sedenions

Robert P. C. de Marrais

**Abstract** Methods for studying zero-divisors (ZD's) in $2^n$-ions generated by Cayley-Dickson process beyond the Sedenions are explored. Prior work showed a ZD system *in* the Sedenions, based on 7 octahedral lattices ("Box-Kites"), whose 6 vertices collect and partition the "42 Assessors" (pairs of diagonals in planes spanned by pure imaginaries, one a pure Octonion, hence of subscript < 8, the other a Sedenion of subscript > 8 and not the XOR *with* 8 of the chosen Octonion). Potential connections to fundamental objects in physics (e.g., the curvature tensor and pair creation) are suggested. Structures found in the 32-ions ("Pathions") are elicited next. Harmonics of Box-Kites, called here "Kite-Chain Middens," are shown to extend indefinitely into higher forms of $2^n$-ions. All non-Midden-collected ZD diagonals in the Pathions, meanwhile, are seen belonging to a set of 15 "emanation tables," dubbed "sand mandalas." Showcasing the workings of the DMZ's (*dyads making zero*) among the products of each of their 14 Assessors with each other, they house 168 fillable cells each (the number of elements in the simple group PSL(2,7) governing Octonion multiplication). 7 of these emanation tables, whose "inner XOR" of their axis-pairs' indices exceed 24, indicate modes of collapsing from higher to lower $2^n$-ion forms, as they can be "folded up" in a 1-to-1 manner onto the 7 Sedenion Box-Kites. These same 7 also display surprising patterns of DMZ sparsity (with but 72 of 168 available cells filled), with the animation-like sequencing obtaining between these 7 "still-shots" indicating an entry-point for cellular-automata-like thinking into the foundations of number theory.

**{0} Review of Prior Results, Motivations for Future Work:** Zero divisors (ZD's), like the "monsters" of analysis before Mandelbrot tamed them into fractals, have been largely avoided as pathological by number theorists. In the last handful of years, however, there has been a shift (largely driven by the needs of physics) toward taking them more seriously. Guillermo Moreno, in 1997, pushed Lie algebraic tactics to the limit of their usefulness and discerned that the ZD's in the 16-D extension, by Cayley-Dickson process, of the Octonions, conform to a pattern homomorphic to the Octonion's automorphism group (and the derivation group of all higher $2^n$-ions from their $2^{n-1}$-ion precursors), $G_2$.[1] His approach, however, was that of an "armchair theorist," in that only a single instance of an actual pair of mutually zero-dividing numbers in the Sedenions was given. K. and Mari Imaeda, meanwhile, focused on manipulating arbitrary (hence complicated) Sedenion expressions, uncovering the ZD's within them, eliciting their general function theory, etc.[2]

Simultaneous with the Imaedas' "top-down" tactics, R. de Marrais took a complementary "bottom-up" approach, isolating underlying structures from which all complicated ZD expressions and spaces in the Sedenions must be composed – structures which "fly under the radar" of classical techniques, requiring an initial calculation-heavy "dirtying of hands" to be elicited and studied.[3] A set of 3 simple "production rules" – and one even simpler generic procedure – were found to suffice for creation of all basic patterns. All points on the pairs of diagonal line elements of 3½ dozen planes (dubbed the "42 Assessors") in Sedenion space mutually zero-divide all other such points on certain other such diagonals. Using the familiar convention of "XOR indexing" (the index of the unit of each of the 15 imaginary axes would equal the XOR of the indices of any 2 units whose product was equal to it), any combination of a pure Octonion (index < 8), and a pure Sedenion (index > 8) which was not the XOR of the Octonion *with* 8, yielded one of 7 x 6 = 42 possible axis-pairs to span Assessor planes.

Representing each such plane by a unique vertex, on one of 7 isomorphic octahedral lattices, resulted in a set of "box-kites" whose 8 triangular faces represented either 4 "sails" sharing vertices with each other but no edges, or 4 empty "vents" in the remaining 4 faces. Tracing an edge along a "sail" connected two mutually zero-dividing diagonals in different Assessors (with points on different diagonals belonging to the *same* Assessor *never* zero-dividing each other). The orthogonal diagonals "/" and "\" of an Assessor's "X," spanned, say, by imaginary axes $i_A$ and $i_B$, would be shorthanded (A+B) and (A-B) respectively, with the Assessor pair as such simply indicated as (A, B).

The multiplication indicated by a box-kite edge would connect two mutually zero-dividing numbers of either same ("/ /" or "\ \") or different ("/ \" or "\ /") orientation, respectively indicated by a "+" or "–" drawn on the joining path. By continuing to trace along a sail's edges, a sequence of mutually zero-dividing pairs (henceforth, "DMZ's," for *dyads making zero*) would be indicated. By tracing all edges of a sail twice over, in the manner of a spinor's "double covering," all 6 diagonals (each making zero with one each from the other two Assessors) would be involved in a single 6-cyclic sequence of DMZ's, in one of two varieties: the 4 sails of any box-kite would partition into 1 "triple zigzag" ("/ \ / \ / \": i.e., "–" signs on all edges; also the case for the opposite-face "vent") and 3 "trefoils" ("/ / / \ \ \": i.e., "+"-signs on all but 1 edge), leading to a total of 6 edges of each sign per box-kite. Other such patterns connecting DMZ sequences, both within and between box-kites, were elicited, and rules were determined for constructing higher-dimensional spaces, all of whose points would zero-divide those of other such spaces.

One more simple property of box-kites must be indicated in this overview: each of the 7 is uniquely associated with the one Octonion which does *not* appear among its 6 vertices. The vertex pairs on opposite ends of any of the 3 orthogonal "struts" stabilizing the framework (like the wooden sticks of a "real-world" box-kite) do *not* form DMZ's – and these are the *only* non-DMZ pairings in a box-kite. Moreover, the products of all 3 of these opposite pairs, depending on the "edge-signs" we might imagine on their struts, display either a pair of the units of index 8 which, per Cayley-Dickson process, *generate* the Sedenions from the Octonions; or, a pair whose index is the "missing Octonion" – henceforth dubbed the given box-kite's "*strut constant*." And, as is easily checked, the "inner XOR" of the Octonion and Sedenion pairs defining the Assessors at each box-kite vertex will always equal '8' plus this strut constant.

In the table below, the strut constant for each box-kite is listed in the first cell of each row, with each vertex's Assessor index-pair listed in the 6 columns that follow. The Octonion indices of the vertices A, B, C are in bold-face: their Assessors form a box-kite's sole triple-zigzag (orange in the diagram, with trefoils red, blue, green). They also form an associative triplet, in "NATO" format (for *natural ascending triplet order*, meaning written so that AB = +C, BC = +A, CA = +B). In fact, each sail's Octonions, read aright, form NATO triplets, viz.: A, D, E; F, D, B; and F, C, E. (Conversely, the 4 triplets contained in any set of 6 Octonions each imply a sail in the box-kite whose strut constant is the 7th.) Note, too, that in each plane spanned by 2 of the 3 orthogonal strut pairs (A, F; B, E; C,D), a circuit of its square frame yields a sequence of 4 DMZ's; but, due to the alternation of signs on the edges, it takes *two* such non-intersecting trips, each linking with a different diagonal in each Assessor, to tour all the square's 8 ZD's.

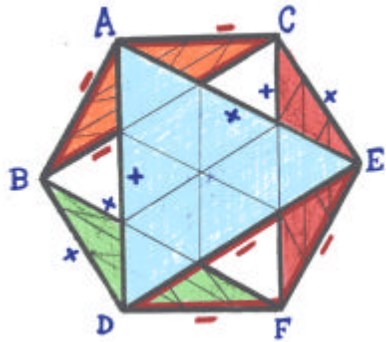

| Strut Const | Assessors at Box-Kite Vertices | | | | | |
|---|---|---|---|---|---|---|
| | A | B | C | D | E | F |
| 1 | **3**, 10 | **6**, 15 | **5**, 12 | 4, 13 | 7, 14 | 2, 11 |
| 2 | **1**, 11 | **7**, 13 | **6**, 12 | 4, 14 | 5, 15 | 3,  9 |
| 3 | **2**,  9 | **5**, 14 | **7**, 12 | 4, 15 | 6, 13 | 1, 10 |
| 4 | **1**, 13 | **2**, 14 | **3**, 15 | 7, 11 | 6, 10 | 5,  9 |
| 5 | **2**, 15 | **4**,  9 | **6**, 11 | 3, 14 | 1, 12 | 7, 10 |
| 6 | **3**, 13 | **4**, 10 | **7**,  9 | 1, 15 | 2, 12 | 5, 11 |
| 7 | **1**, 14 | **4**, 11 | **5**, 10 | 2, 13 | 3, 12 | 6,  9 |

Before offering more pictures and describing further properties, let's consider just how peculiar (and hence, perhaps useful) the results just described truly are. Specifically, let's see how radically they differ from the only zero-dividing entities familiar (indeed, indispensable) in fundamental physics. Though their couching in operator or Clifford algebra formalism may hide their ultimate nature as *number types* from casual view, the "projection" operators of quantum mechanics are, in fact, hybrid units in a *toroidal* number space. (Imaginary, Quaternion, Octonion, and higher $2^n$-ion units, by contrast, run through their paces in 2-, 4-, 8-, and $2^n$- dimensional *spherical* spaces.)

The simplest such number-form torus is just the Cartesian product of two imaginary orbits: the study of the 4-D result of two 2-D spaces containing the unit circles of the Argand diagram. The Italian geometer Corrado Segrè first studied these "bicomplex" numbers in the 1890's[4]; G. Baley Price wrote the first full-length text studying them as such, and in depth, only in the *19*90's.[5] (They also are implicit in Roger Penrose's well-known "twistor" formalism[6], and were the starting point of Charles Musès' less well-known "hypernumber" program.[7]) Segrè noted the 4 axes defining their context were generated by 2 imaginary units, the real unit, and a fourth "mirror number" whose unit was a square root of +1 but not itself equal to $\pm 1$. This last is well-known in its own right: a space-spanning system of 3 such, mutually noncommuting in the manner of the Quaternions from which one can derive them, are commonly called the "Pauli spin-matrices." And as Musès' work (and the detail-work of the 8-D Clifford algebra $Cl_3$[8]) underscores, a "normal" imaginary (or quaternion) unit doesn't commute with any such "mirror units" – but Segrè's second imaginary (represented by the 2 x 2 identity matrix, but with the usual **i** replacing each real unit on the main diagonal) commutes with both sorts.

If we call the mirror unit **m**, the two diagonals in the (1,**m**) plane each have a point which is *idempotent*: $[½(1\pm\mathbf{m})]^2 = ¼ (1 \pm 2\mathbf{m} + 1) = ½ (1\pm\mathbf{m})$; moreover, the bracketed quantities, raised to *arbitrary* powers, remain unchanged, allowing for exponential and other functions to be defined along the diagonals which have these points as their "units" in a manner isomorphic to the workings of the real number line. These points also have two other properties of note: 1) they are the centers of mutually orthogonal Argand-diagram-like "power orbit" circles, with diameters extending from $\pm\mathbf{m}$ to +1, allowing them to provide a second "S1 x S1" representation of the generation of the torus containing them; and, 2) the *product* of these mutually orthogonal idempotents, as is readily checked, is *zero* – making them *nil*potent as well. In this latter aspect, they are easily recognized as the number-form basis for "projection" operators.

Such toroidal numbers, far from being "artificial," are forced upon us by nature's penchant for balanced bookkeeping: they're what you get when you contemplate the phase curves defined by linear equations with purely imaginary eigenvalues (e.g., in the theory of oscillations of conservative systems). But perhaps it's artificial to consider such projections into the Quantum Void as necessarily built up from such zero-divisor dyads. Perhaps this restriction has something to do with the scandal of the universe's "90% missing mass"? Spherical number forms, beginning with the Sedenions, have ZD's on diagonals, too – but they're never idempotents, have no real components, obey a ternary not binary logic, and those composed from the *same* units will *never* (as opposed to *always*) form DMZ's. The difference between these two approaches is as radical as that between the renormalization tactics used in plus-or-minus-charged, photon-gauged, electromagnetism, and the three-color "strong force" of quarks and gluons.

Among other possibilities, this just indicated difference has a bearing on how one interprets the well-known "zero-point fluctuations" (ZPF) of the vacuum – as manifest in the experimentally well-attested "Casimir force," which allows actual extraction of energy (admittedly in minuscule quantities!) from the Void. And this, in turn, soon involves one in foundation questions concerning the interrelations of the quantum and general relativistic realms. Standard quantum field theory represents ZPF via standard creation and annihilation operators acting on the vacuum; and, as part of the same gesture, gravity is treated as an exchange of so-called "gravitons" in *flat* [sic] spacetime.

Alternative approaches are numerous; we note but two, which point back independently to a 1968 "blue sky" piece by Andrei Sakharov[9]. First is the deep investigation initiated by Haisch, Rueda and Putoff,[10] of a non-Machian, zero-point basis for *inertia*, developed using the classical, Poynting-vector-based tools of stochastic electrodynamics. In their view, mass can be seen as a way of characterizing the resistance due to the ZPF (and general Quantum Vacuum) "molasses" that kicks in upon acceleration. (The asymmetric distortion induced in the ZPF by accelerating bodies was demonstrated, as they note, by Davies and Unruh in the mid-'70s, building on Hawking's "black-hole evaporation" model.) They hypothesize that rest mass would similarly be associated with the fundamental "restlessness" of a particle (Schrödinger's "zwitterbewegung"), attaining to resonance at the Compton frequency, thereby finding the origin of de Broglie's "matter waves" in Doppler shifts of the ZPF spectrum. They further suggest expanding their approach to peel back, in a manner first suggested by Sakharov, the formalism of general relativity itself. Which segués into the second alternative approach.

In a book-length investigation of <u>Nonassociative Algebras in Physics</u>, Jaak Lõhmus, Eugene Paal and Leo Sorgsepp note two patterns they assume are interconnected.[11] On the one hand, they observe that modern physics has built itself upon the recognition of a sequence of "infobarriers" associated with ever less regular notions of number: with relativity, imaginaries are associated with the restrictive causality of the light-cone; the Heisenberg uncertainty relations depend upon the Quaternions' noncommutativity; arguments are made for such phenomena as "exact color symmetry" and its concomitant quark confinement being side-effects of the Octonions' nonassociativity.

On the other hand, they note that attempts to quantize gravity largely fail due to the unique nature of the gravitational force: the analog of the "fine-structure *constant*" is in fact a *variable*. Phenomena at the ultra-small scale of the Planck length require a second kind of quantization – not of action, but of *distance itself*, leading to a notion of a spacetime "crystalline" structure. But it is almost tautological that this latticework is volatile: below the Planck length, it decays, and "multiplicands" are somehow "elastic."

Putting these two threads together, it seems hardly outrageous to suggest that Sakharov's "metrical elasticity" underwriting gravity is directly connected to the *next* breakdown of properties – the loss of a field-theoretic notion of *norm* by which ZD's are spawned. Since *all* points of a ZD diagonal make DMZ's with *any* points of any other suitably box-kite-connected ZD, can we not say a study focused on the "pre-metric" latticework of ZD's, in the Sedenions and beyond, provides the *only* logical way to explore such phenomena – and, perhaps identically, the phenomenology of the "Creation pressure" leading up to the Big Bang? The connection, after all, of the Sedenions to the Dirac equation, and hence quantum mechanics, is well known; and a link from the Sedenions' ZD substrate to general relativity can also be strongly suggested, thanks to box-kites.

The following diagram of a box-kite doesn't have signed edges, and has peculiar symbols on its vertices, and its "sails" lack coloring; it is otherwise, however, indistinguishable from that provided above. This was hardly the intent of its creator, though: for it is an aid in visualizing a key proof concerning the symmetries of the *curvature tensor* in John Milnor's classic lecture notes on <u>Morse Theory</u>.[12] Using the notation of affine connections, Milnor states in a lemma that a Riemannian manifold's curvature tensor satisfies four relations: (1) it is skew-symmetric ($R(X,Y)Z + R(Y,X)Z = 0$); (2) cyclically ordered such triples of vector fields sum to zero ($R(X,Y)Z + R(Y,Z)X + R(Z,X)Y = 0$); (3) bracket relations are skew-symmetric ($<R(X,Y)Z,W> + <R(X,Y)W,Z> = 0$) in outer symbols, but (4) unchanged when inner and outer symbols are exchanged ($<R(X,Y)Z,W> = <R(Z,W)X,Y>$). The fourth property is the most difficult, its proof requiring use of the other three – and leading to the diagram as an aid to its following. Says Milnor:

> Formula (2) asserts that the sum of the quantities at the vertices of shaded triangle W is zero. Similarly (making use of (1) and (3)) the sum of the vertices of each of the other shaded triangles is zero. Adding these identities for the top two shaded triangles, and subtracting the identities for the bottom ones, this means that twice the top vertex minus twice the bottom vertex is zero. This proves (4), and completes the proof.

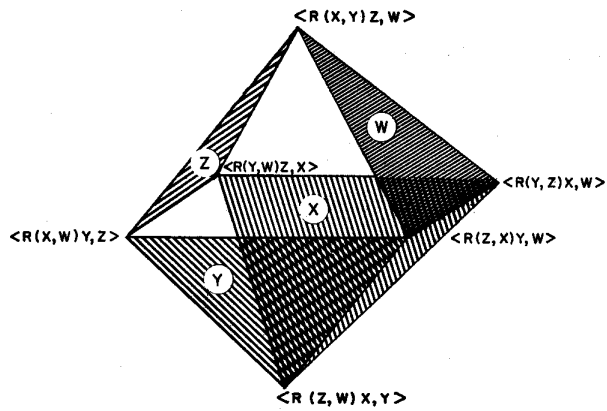

With Sedenion-based box-kites, of course, the shaded triangles comprise the sails, and the *products* along each of their edges are zero; the terms at opposed ends of struts, meanwhile, obey a different relation from those connected to each other in a sail – not their sums, but their *differences*, must be attended to. Can a small number of "symmetry breaks" and/or "categorical removes" serve to link these radically different uses of the box-kite diagram? The question is a deep one, and will be treated in follow-up studies. For now, consider that: (1) tensor notation replaced the now largely extinct approach of Dyads and Nonions; (2) J. J. Sylvester's Nonions focused on a number notion he thought, from the vantage of determinant theory, exactly analogous to Quaternions, but involving *cube* roots of $+1$ rather than *square* roots of $-1$, entailing 3 x 3 arrays of monomials associated with 3 x 3 matrix representatives whose non-zero entries were 1's and imaginary cube roots of unity[13]; (3) the monomials are precisely those which "Siersma's Trick" in Catastrophe Theory displays for the Double Cusp[14]; (4) ZPF "drag force" is a *third order* effect; and (5) while un-normed Sedenions do not obey an "n-squares law" like the Imaginaries' (2, 2, 2), they permit, among other relations of interest, a (9, 16, 16) variant.[15]

Only the general context of that fifth and last consideration will be dealt with herein. For, properly appreciated, it suggests a "metalevel" from which to view quantum phenomena comparable to that provided for viewing general relativity and the relation of mass to inertia in the above discussion. The problem of "n-squares laws" generalizes the Pythagorean theorem and its application to higher-dimensional analogs of *distance*: the (2, 2, 2) rule for Imaginaries merely says, firstly, that any Imaginary can be construed as a point on a circle, hence at a distance of a "radius" from the origin, so that the radius' square is the sum of the squares of the abscissa and ordinate values of the same point; and, secondly, that the product of *two* such points can *itself* be understood the same way.

More elaborate extensions of the same Pythagorean logic lead to the (4, 4, 4) rule for Quaternions (the product of one sum of four squares times another is yet another sum of four squares), and, more elaborate still, to the (8, 8, 8) rule for Octonions. Hurwitz's classic proof, though, of the breakdown of the norm (and hence field structure) beyond 8 dimensions, means that Sedenions are most readily assimilated to a (16, 16, 32) rule: since any Sedenion can be partitioned into two sets of 8 coordinates A and B, the product of two can be written $(A+B) \cdot (A'+B') = AA' + AB' + BA' + BB' = 4 \times 8 = 32$ dimensions required. And as Cayley's friend, the Rev. T. Kirkman, showed more than a century and a half ago, products of *three or more* Sedenions actually require a backdrop of *64 dimensions*![16]  All of which suggests that the indefinite cascades of virtual particles, handled by renormalization methods, may perhaps signify something else entirely, but "seen through a glass darkly": viz., the 8D-chunked (or, "going toroidal" in E8 x E8 superstring style, 8D + 8D-chunked) *side effect* of a "Bizarro World" mutation of Hilbert space thinking, where a proliferation of ever more multiply connected *ZD-nets* of $2^n$-ions, n approaching infinity, provides (at least the "cover story" for) the "hidden agenda."

We will return to the general architectural issues implied here at the end of this monograph, when we consider the Artificial-Life- or neural-net- like simulation-work suggested by such blue-sky "Creation Pressure" thematics. What is clear here, though, is that even broaching such topics makes adherence to the tried and true (and badly aging) Lie algebraic toolkit seem (if not exactly *laughable*) inadequate. The best way to make this patently obvious is to thoroughly *and concretely* investigate the realities being only vaguely indicated by the highly useful (but clearly limited) "$G_2$" results of Moreno. In later work than that cited earlier, he determines that the automorphism group of the ZD's of all $2^n$-ions derived from the Cayley-Dickson process, from the Sedenions on up, obey a simple pattern: for $n \geq 4$, this group has the form $G_2 \times (n-3) \times S3$ (this last being the standard order-6 permutation group on 3 elements).[17]

This says the automorphism group of the Sedenions' ZD's has order $14 \times 1 \times 6 = 84$: the number of diagonal line elements, all of whose points form DMZ's only with other such line elements, contained in the 42 Assessors. Let's expand our horizons to the 32-D and 64-D cases and give them names (since they currently don't have any!). Following the convention adopted by Tony Smith, who has called the $2^8$-ions **Voudons** after the 256 deities of the Ifa pantheon of Voodoo or Voudon, we'll call these **Pathions** (for the "32 Paths" of Kabbalah) and **Chingons** (for the 64 Hexagrams of the I Ching or Book of Changes). For completeness' sake, we name the $2^7$-ions too, dubbing them **Routons**, after that legendary source of high-tech innovativeness, *Route 128* of the "Massachusetts Miracle" that paralleled Silicon Valley's on the "Left Coast" of this country. From Reals to Voudons, then, we extend the standard shorthand thus: **R**, **C**, **H**, **O**, **S**, **P**, **X**, **U**, **V**.

**{1} Using "Table-Driven" Cayley-Dickson to Excavate the Kite-Chain Midden:**

What is meant by "table-driven": in the early days of the PC Revolution (the one focused on computing, not "correctness"), dBase II allowed you to do easy database programming, but largely confined you to handling data at the *record* level. When the next wave of tools emerged, Paradox (the precursor of the networked, SQL-based juggernaut that Oracle launched) let you use a simple logic of columnar links and commands to operate on entire *tables* of data. In the transition phase between technologies, it was not uncommon to take dBase programs that took an hour to run and replace them with Paradox routines that did the same job in under a minute. The moral of this story is: the way the Cayley-Dickson process (henceforth, CDP) is typically presented (and, apparently, used) is like dBase II. If we're to do anything interesting with ZD patterns in $2^n$-ions, n large, we must follow the example of Paradox. (And, of course, we must also drop our pencils at some not so far off point and write programs!) We will soon see that the "table-driven" approach to the CDP is ridiculously easy to employ and explain.

What the "Kite-Chain" (pronounced "kitchen") "Midden" is: the PC Revolution was roughly contemporary with the spread of the "object-oriented programming," or "OOP" paradigm. Its way of viewing things provides a handy (and deep-running) metaphor for thinking about zero divisors. Two key notions in OOP lore are design-hiding and garbage-collecting (the latter conveniently built into the Java language). The first embodies the precept that a maintainable and reusable design module will "privatize" key variables, denying external users direct access to them. Beginning at least with "quark confinement," Nature seems to make its innermost secrets likewise difficult of access. However, Nature – from the panda's thumb to the microwave background radiation of the cosmos – also does an imperfect job at covering up its tracks: pragmatic trade-offs being what they are, allocation and deallocation of memory is never perfectly handled by any "garbage collection" algorithm, in Nature or in software. Indeed, one way of looking at archaeology (and in this sense, cosmologists studying the faint traces of the Big Bang's aftershocks are engaged in it) is that this imperfection makes it possible.

A *midden*, in Old English, comes from the word for "dung heap," and has come to mean, to archaeologists, a specific thing: the heaps of broken dishes, worn-out cutlery, sabre-tooth tiger bones, and other refuse tossed into deep recesses of inhabited caves or other primitive dwelling sites over the millennia: such *kitchen middens* provide a treasure trove of information on our ancestors' digestive tracts, sophistication of toolcraft, and countless other matters that students of the past crave to (literally) unearth.

Divisors of zero are likewise a "dung heap" of number-theory leavings which no one's had much use for; yet it is possible to excavate their layers and find some surprising patterns which extend indefinitely into remote reaches not of antiquity, but of "$2^n$-ionity" (which indeed might signify the ultimate antiquity of things leading up to and immediately out of the Big Bang anyway). The endless succession of layers in this "dig" takes the form of a never-completed chain of box-kite-like entities, hence the "kite-chain" (said "kitchen") qualifier on the "midden." With a table-driven CDP in hand, the intricacies of the *kite-chain midden* will be easy to elicit. And once we've grasped their workings, the patterns they reveal to us will make us doubt the generality of the "$G_2$" approach! (Especially when we see, in the section after this one, that as badly as Moreno's "upper bound" is exceeded by the kite-chain midden, a slew of more elaborate ZD structures we'll call "sand mandalas," first appearing in the Pathions, exceed such limits even more!)

How the CDP works: like the route to Chaos we think the CDP resembles, there is a "period-doubling" (here, doubling of *dimension*) which the CDP algorithm captures by partitioning the numeric entities of its next level up in two parts: units whose indices are contained in the starting-point $2^n$-ions go to the left of the comma, and those whose indices range from $2^n$ to $2^{n+1}-1$, go to its right – but with $2^n$ subtracted from them. The complex numbers, to start (as the CDP does) with the simplest case, have n = 0, and are inscribed (a,b), with a and b scalars multiplying the units $\mathbf{i_0}$ and $\mathbf{i_1}$ – reals and imaginaries.

The product of two complex numbers, written (a,b) and (c,d), expands and recontracts as follows: $(a\mathbf{i_0} + b\mathbf{i_1})(c\mathbf{i_0} + d\mathbf{i_1}) = (ac - bd)\mathbf{i_0} + (ad + bc)\mathbf{i_1} = (ac - bd, ad + bc)$. CDP doubling yields the Quaternions, which obey much the same pattern. But since their units don't commute, order matters. And since there's more than one imaginary unit to deal with, indices must be made explicit, and so must the higher-dimensional analog of complex conjugates. With simple imaginaries, the conjugate of $\mathbf{C} = (a\mathbf{i_0} + b\mathbf{i})$ is just $\mathbf{C^*} = (a\mathbf{i_0} - b\mathbf{i})$, so that $\mathbf{C \cdot C^*} = a^2 + b^2$; for general Quaternion $\mathbf{Q} = (a\mathbf{i_0} + b\mathbf{i_1} + c\mathbf{i_2} + d\mathbf{i_3})$, the conjugate $\mathbf{Q^*}$ is likewise $(a\mathbf{i_0} - b\mathbf{i_1} - c\mathbf{i_2} - d\mathbf{i_3})$, and $\mathbf{Q \cdot Q^*} = a^2 + b^2 + c^2 + d^2$. In CDP notation, we get $(a + b\mathbf{i_1}, c + d\mathbf{i_1})$ and $(a - b\mathbf{i_1}, -c - d\mathbf{i_1})$ respectively. For general $2^n$-ions, the $\mathbf{M}$ and $\mathbf{M^*}$ always have equal components, of opposite sign for non-real terms only.

Here, we see that the $\mathbf{i_2}$ term (and, more generally, the term whose index is the next-higher power of 2) acts as the *generator* of the new units appearing at the next level (hence, has the same index – 0 – on the *right* of the comma that the real unit has on the *left*); as such, a generator *always* produces same-signed resultants when applied on the same side of an arbitrary unit among the $2^n$-ions it generates. Given that one can't index even the Octonions so that all NATO triplets are also in counting order – (1, 7, 6) and (3, 6, 5) are out of order in the usual XOR indexing scheme, for instance – this implies that no imaginary unit *except* the generator will always produce a same-signed resultant with arbitrary other units multiplied from the same side. This special feature of generators has far from trivial ramifications with higher $2^n$-ions, as we'll see soon enough.

We need note one more general feature before giving the standard iterative CDP algorithm, for it's more general than we need. We can generate higher-order systems of mirror, imaginary, or "dual" numbers (which square to 0), by setting a parameter to +1 or – 1 or 0 respectively, and create mixes of such systems by varying the parameter at each successive CDP iteration. (One can even set the last of 4 parameters to a radical, to embed the Sedenions in an extension field which has no zero divisors at all. More generally still, one can exploit the theory of so-called "Pfister forms" – or engage in category theoretics to similar effect – and define the CDP over polynomial fields, expunging zero-divisors and tying $2^n$-ions to quantum groups in one fell swoop[18]: such nifty tricks, however, won't concern us.) For our purposes, all parameters will be  –1's; for completeness, we give the general expression, with the most recently iterated parameter written **p**:

$$X \cdot Y = (A,B) \cdot (C,D) = (AC + \mathbf{p}(D^*)B, B(C^*) + DA)$$

If the component parts are *not* complicated admixtures, but simple sets of non-0-indexed imaginary units with unit coefficients, we can treat conjugation as just negation, and write only indices and signs (M – N, instead $a\mathbf{i_M} + b^*\mathbf{i_N}$). Which gives us this:

$$X \cdot Y = (A,B) \cdot (C,D) = (AC + DB, DA - BC)$$

Before we go into "table-drive," let's make this real with a worked example. Going back to the Box-Kite diagram and its associated table, let's take a random path between vertices – BC, say, on the triple zigzag ABC. Now, let's pick a random strut constant: 6, for instance. The table says B then equates to the Assessor (4,10), and C to (7,9). The sign on the edge is minus, so let's flip a coin and pick the "+" diagonal to be B, the "–" to be C: this means multiply ($i_4 + i_{10}$) by ($i_7 - i_9$). By the formula, we get:

$$(4, 2) \cdot (7, -1) = (\,[(4)(7) - (1)(2)],\ [-(1)(4) - (2)(7)]\,)$$

Refer to the relevant "NATO triplets": (3, 4, 7); (1, 2, 3); (1, 4, 5); (2, 5, 7). The CDP comma separates **O**'s from **S**'s, the latter generated from the former by the unit indexed $2^3 = 8$, so we get $(3 - 3, -5 - (-5)) = (+i_3 - i_{13}) + (- i_3 + i_{13})$. For "hard" parentheses, "pair creation from the Quantum Void" is suggested; for "soft," we just have **0**.

Note that 3 *xor* 13 = 6, the strut constant, as required. Furthermore, (3, 13) has edge sign "–", so we can continue edge-multiplying and get further zeros in a DMZ 6-cycle, and say that the Assessor (3, 13) at vertex A is "emanated" by Assessors (4, 10) and (7, 9), with $(4 + 10) \cdot (3 - 13) = (3 - 13) \cdot (7 + 9) =$ [pair creation or] 0, and so on.

It's easy to find multiplication tables for the Octonions, and also easy to build one for the Sedenions; no higher tables have been published that we know of, but we don't need them. All we need do is realize that the table for $2^{n+1}$-ions is comprised of 4 tables of the same size as that for the $2^n$-ions, with the latter table constituting the upper left quadrant of the larger ensemble, its row and column values running from 0 to $2^n - 1$. Calling this quadrant "*I*," and continuing to designate the other 3 by successive Roman numerals in clockwise order, we can fill in all values of each of these by a simple transformation of the ensemble of $2^n$-ion values in *I*.

*Quadrant II:* Row values same as *I*; column values of $n^{th}$ column of *II* = those for $n^{th}$ column of $I + 2^n$. Hence if product of row-unit *r* and column-unit *c* in *I* is written $(r, 0)(c, 0) = (rc + 0, 0 \cdot r - 0 \cdot c) = rc$ = entry in cell (r, c), then cell (r, $2^n$ + c) in *II* will be $(r, 0)(0, c) = (r \cdot 0 + c \cdot 0, c \cdot r - 0) = (0, -rc) = -(2^n + rc)$, with "special handling" for *r* or *c* = 0 or *r* = *c*. In the last case, since any $i_n^2$ is negative unity, the cell entry will be the odd-looking "–0" (short for $-i_0 = -1$) in the long diagonal of *I*. The peculiar status of the *generator* – any smaller-indexed positive unit *u* multiplying it on the left yields the *positively* signed unit whose index is $2^n + u$ – means that the corresponding long diagonal in *II* will show $(-2^n)$ for all entries save the first row's. Corollarily, all entries in the first row and column will be *positive*, with content $2^n$ greater than the corresponding cell in *I*.

To simplify, collectively refer to the extension of the top row, the first new ($2^{nth}$) column, and the main diagonal of *II* the "trident" of *II* (with tridents in all quadrants defined analogously). Then we can sum up this way: the *index-values* in cells (r, $2^n$+c) = $2^n$ + (r, c); *sign* of cells (r, $2^n$+c) = *same* as for (r, c) if in trident, *opposite* otherwise.

*Quadrant III:* by similar arguments, cell entries in trident of *III* are *identical* in content to corresponding entries in trident of *I*, and *same* in sign *except for the leftmost* ($2^{nth}$) *column*, where signs are *reversed*; for all other cells, ($2^n$ + r, $2^n$ + c) = *same* in content and *opposite* in sign as (r, c).

*Quadrant IV:* for each entry ($2^n$ + r, c) content is *identical* to that in corresponding cell (r, $2^n$ + c) of *II*, with signs *reversed* in the trident (*except* for the extension of $0^{th}$ column), but the *same* everywhere else.

Naturally, the best way to grasp the iterative workings of this process is to examine a concrete for-instance where it's been applied multiple times. On the next page, we give the multiplication table for the Sedenions (the Pathions' 32 x 32 array won't fit), using a checkerboard scheme of shading alternate Quaternion-sized blocks of cell entries, to make the iterative structure easier to parse visually. Beneath this table, we give a listing of the $(2^5-1)(2^5-2)/3! = 155$ NATO triplets in the Pathions (which we'll be focused on most), with the triplets first appearing in the Octonions and Sedenions broken out separately at the top of the listing. Beyond the Pathions, even a listing of triplets becomes unmanageable, as there are $(2^6-1)(2^6-2)/3! = 651$ of them for the Chingons alone, so readers interested in checking calculations involving them later in this argument will need to feed Pathion triplets into the CDP formula (or write code!).

One crucial matter the "table-driven" approach makes obvious: since no box-kite vertices have "trident" entries, one gets box-kite "harmonics" simply by adding 16 (or 32, or 48 … ) to their indices. These entities are not "stand-alone," however: for o, o', o" Octonions; S, S', S" Sedenions; n an integer; two "Assessor harmonics" with indices (o + 16n, S + 16n) and (o' + 16n, S' + 16n) will emanate some (o" [ = o *xor* o' = S *xor* S' ], S" [ = o *xor* S' = S *xor* o' ] ) contained in the basic box-kite. Each of the 12 paths in such a harmonic, then, will have its "sail" completed by the Assessor in the "base-line" box-kite whose harmonic, by visual logic, "ought to" complete its sail, but doesn't.

Stranger still, if we consider all 3 harmonics allowed when we double domains up to the 64-D Chingons, not only will the base-line box-kite complete sails in each, but the first harmonic will form an "overtone series" with the second and third: "sails," that is, will have one "Assessor harmonic" from each! For example, if the base-line is the box-kite with strut constant 3, its Assessor strut-pairs are (2, 9) and (1, 10), (5, 14) and (6, 13), and (7, 12) and (4, 15) – the first-listed in these pairings forming the box-kite's triple-zigzag. The "Assessor harmonics" for n = 1, 2, 3 respectively are therefore, in same order:

| | | |
|---|---|---|
| (18, 25) and (17, 26) | (21, 30) and (22, 29) | (23, 28) and (20, 31) |
| (34, 41) and (33, 42) | (37, 46) and (38, 45) | (39, 44) and (36, 47) |
| (50, 57) and (49, 58) | (53, 62) and (54, 61) | (55, 60) and (52, 63) |

"Sails" involving all 3 of these harmonics can be constructed in sets of four, each set completed by one of the 6 Assessors in the first harmonic. For the lowest-index instance, the four sails are these:

| | |
|---|---|
| (17, 26) – (36, 47) – (53, 62) | (17, 26) – (37, 47) – (52, 63) |
| (17, 26) – (38, 45) – (55, 60) | (17, 26) – (39, 44) – (54, 61) |

If we only consider the first harmonic – all we get with just the Pathions – we have a system of twice the number of ZD-diagonals in a base-line box-kite, all interrelated in the manner just described. This fits quite nicely with Moreno's assessment of the total count of "irreducible" zero-divisors we should have: per his formula, $2^5$-ions should have (5-3) x 84 = 168 such entities. But clearly, the formula breaks down with one more doubling: when the kite-chain midden is extended to $2^6$-ions, we get not 3 x 84, but 4 x 84 = 336 ZD-diagonals: the maximal symmetry of Felix Klein's famous figure that Poincaré determined was a tessellation by triangles of the hyperbolic plane.[19]

|    | 0  | 1   | 2   | 3   | 4   | 5   | 6   | 7   | 8   | 9   | 10  | 11  | 12  | 13  | 14  | 15  |
|----|----|-----|-----|-----|-----|-----|-----|-----|-----|-----|-----|-----|-----|-----|-----|-----|
| 0  | 0  | 1   | 2   | 3   | 4   | 5   | 6   | 7   | 8   | 9   | 10  | 11  | 12  | 13  | 14  | 15  |
| 1  | 1  | -0  | 3   | -2  | 5   | -4  | -7  | 6   | 9   | -8  | -11 | 10  | -13 | 12  | 15  | -14 |
| 2  | 2  | -3  | -0  | 1   | 6   | 7   | -4  | -5  | 10  | 11  | -8  | -9  | -14 | -15 | 12  | 13  |
| 3  | 3  | 2   | -1  | -0  | 7   | -6  | 5   | -4  | 11  | -10 | 9   | -8  | -15 | 14  | -13 | 12  |
| 4  | 4  | -5  | -6  | -7  | -0  | 1   | 2   | 3   | 12  | 13  | 14  | 15  | -8  | -9  | -10 | -11 |
| 5  | 5  | 4   | -7  | 6   | -1  | -0  | -3  | 2   | 13  | -12 | 15  | -14 | 9   | -8  | 11  | -10 |
| 6  | 6  | 7   | 4   | -5  | -2  | 3   | -0  | -1  | 14  | -15 | -12 | 13  | 10  | -11 | -8  | 9   |
| 7  | 7  | -6  | 5   | 4   | -3  | -2  | 1   | -0  | 15  | 14  | -13 | -12 | 11  | 10  | -9  | -8  |
| 8  | 8  | -9  | -10 | -11 | -12 | -13 | -14 | -15 | -0  | 1   | 2   | 3   | 4   | 5   | 6   | 7   |
| 9  | 9  | 8   | -11 | 10  | -13 | 12  | 15  | -14 | -1  | -0  | -3  | 2   | -5  | 4   | 7   | -6  |
| 10 | 10 | 11  | 8   | -9  | -14 | -15 | 12  | 13  | -2  | 3   | -0  | -1  | -6  | -7  | 4   | 5   |
| 11 | 11 | -10 | 9   | 8   | -15 | 14  | -13 | 12  | -3  | -2  | 1   | -0  | -7  | 6   | -5  | 4   |
| 12 | 12 | 13  | 14  | 15  | 8   | -9  | -10 | -11 | -4  | 5   | 6   | 7   | -0  | -1  | -2  | -3  |
| 13 | 13 | -12 | 15  | -14 | 9   | 8   | 11  | -10 | -5  | -4  | 7   | -6  | 1   | -0  | 3   | -2  |
| 14 | 14 | -15 | -12 | 13  | 10  | -11 | 8   | 9   | -6  | -7  | -4  | 5   | 2   | -3  | -0  | 1   |
| 15 | 15 | 14  | -13 | -12 | 11  | 10  | -9  | 8   | -7  | 6   | -5  | -4  | 3   | 2   | -1  | -0  |

**Octonions** *(1 2 3)* *(1 4 5)* *(1 7 6)* *(2 4 6)* *(2 5 7)* *(3 4 7)* *(3 6 5)*
**Sedenions** *(1 8 9)* (1 11 10) (1 13 12) (1 14 15)
*(2 8 10)* (2 9 11) (2 14 12) (2 15 13) *(3 8 11)* (3 10 9) (3 13 14) (3 15 12)
*(4 8 12)* (4 9 13) (4 10 14) (4 11 15) *(5 8 13)* (5 10 15) (5 12 9) (5 14 11)
*(6 8 14)* (6 11 13) (6 12 10) (6 15 9) *(7 8 15)* (7 9 14) (7 12 11) (7 13 10)
**Pathions**
*( 1 16 17)* (1 19 18) (1 21 20) (1 22 23) (1 25 24) (1 26 27) (1 28 29) (1 31 30)
*( 2 16 18)* (2 17 19) (2 22 20) (2 23 21) (2 26 24) (2 27 25) (2 28 30) (2 29 31)
*( 3 16 19)* (3 18 17) (3 21 22) (3 23 20) (3 25 26) (3 27 24) (3 28 31) (3 30 29)
*( 4 16 20)* (4 17 21) (4 18 22) (4 19 23) (4 28 24) (4 29 25) (4 30 26) (4 31 27)
*( 5 16 21)* (5 18 23) (5 20 17) (5 22 19) (5 25 28) (5 27 30) (5 29 24) (5 31 26)
*( 6 16 22)* (6 19 21) (6 20 18) (6 23 17) (6 25 31) (6 26 28) (6 29 27) (6 30 24)
*( 7 16 23)* (7 17 22) (7 20 19) (7 21 18) (7 26 29) (7 27 28) (7 30 25) (7 31 24)
*( 8 16 24)* (8 17 25) (8 18 26) (8 19 27) (8 20 28) (8 21 29) (8 22 30) (8 23 31)
*( 9 16 25)* (9 18 27) (9 20 29) (9 23 30) (9 24 17) (9 26 19) (9 28 21) (9 31 22)
*(10 16 26)* (10 19 25) (10 20 30) (10 21 31) (10 24 18) (10 27 17) (10 28 22) (10 29 23)
*(11 16 27)* (11 17 26) (11 20 31) (11 22 29) (11 24 19) (11 25 18) (11 28 23) (11 30 21)
*(12 16 28)* (12 21 25) (12 22 26) (12 23 27) (12 24 20) (12 29 17) (12 30 18) (12 31 19)
*(13 16 19)* (13 17 28) (13 19 30) (13 23 26) (13 24 21) (13 25 20) (13 27 22) (13 31 18)
*(14 16 30)* (14 17 31) (14 18 28) (14 21 27) (14 24 22) (14 25 23) (14 26 20) (14 29 19)
*(15 16 31)* (15 18 29) (15 19 28) (15 22 25) (15 24 23) (15 26 21) (15 27 20) (15 30 17)

**{2} The Pathions' "Sand Mandalas" (Surprises in Symmetry and Degeneracy):** In the listing of triplets just given, those none of whose 3 pairings can act as Assessors (i.e., contribute ZD's) are rendered in italics. These include all the Octonions, all Sedenions with an '8,' and all Pathions with a '16.' And, just as each of the 7 x 3 "8-less" Sedenion triplets of form (O, S, S') produced two Assessor index-pairs – (O, S) and (O, S') – making 42 in all, so the Pathions have 15 x 7 = 105 "16-less" triplets of form (O or S, P, P'), and which yield 210 "Assessors" organized in 15 ensembles which are the Pathion equivalent of Box-Kites. In addition, as implied in the above discussion of kite-chain middens, 42 of the (P, P') pairings are ZD's as well, which makes 252 x 2 diagonals in each = 504 "irreducible" ZD's specific to the Pathions. Add in the 84 first showing in the Sedenions, and we have **588** in all: 3½ times the count Moreno's formula would allow.

The 15 Pathion ensembles of Assessors have many features we've already seen in box-kites: Assessors at all 14 vertices of the analogous lattice have the same inner *xor*, its value – one of 15 possible (O or) S indices – acting as the ensemble's signature. This signature plays the role of "strut constant" as well, so that the Assessors at the 14 vertices are arrayed in opposing pairs which do *not* form DMZ's, since their products include non-canceling pairs of strut-constant-indexed and index-16 imaginaries. And, as with box-kites, not just the (here, index- **16** instead of **8**) generator, but the Assessor whose Pathion is the *xor* of the strut constant with it, is also excluded from the lattice.

While this much is easy to see, visualizing the connectedness of these ensembles is hard: with each having 14 vertices, possibly forming DMZ's with all but 1 of the others, we have 84 edges, forming (14 x 12)/3! = 28 "sails," divvied out among 7 Box-Kites, each of which shares each of its 3 strut with 2 of the others. We say "*possibly* forming" advisedly: one of the most surprising features of these ensembles – and why, as we'll see, we call them "sand mandalas" – is that the full count of 168 DMZ's is only obtained when the strut constant is 8 or less. And, when it *is* precisely 8, we get a peculiar phenomenon we'll call "edge polarization." To grasp these things, however, a different sort of visualizing strategy is called for. In a partial generalization of the table appended to the box-kite diagram toward the start of this paper, we'll display the behaviors of the sand mandalas in 14 x 14 spreadsheet-like layouts we'll call *emanation tables*.

As with the column headings there, indices will be arranged to best satisfy two sometimes conflicting constraints. First, Assessors on opposite ends of a strut are given symmetrically opposite placement, so that indices have "nested parentheses" ordering (just as columns A and F, B and E, and C and D formed strut-opposed pairings in the earlier table). Second, when this first constraint allows (as it does for the first half of the column heads for all sand mandalas with strut constants of 8 or greater), the lower-indexed Assessors from all strut-pairs should have their O or S terms placed atop the columns in ascending left-to-right order (and top-to-bottom order on the rows).

When arranged in this manner, all cells of both long diagonals of the table will always be empty: those running from the upper left toward the bottom right represent products of an Assessor's diagonals with themselves, and these can never form DMZ's; those running downward from the upper right, meanwhile, represent products of "strut opposites," and these will always contain pairs of index-16 and strut-constant-indexed imaginaries. The remaining $14^2 - (2 \times 14) = 168$ cell entries are assumed, by convention, to be the O or S indices of the ZD "emanated" by multiplying the row-head's Assessor times the column-head's, with appropriate edge-sign prefixing the entry.

We'll construct the first "emanation table," for strut constant = 1, now. All the Assessors have inner *xor* 17, units indexed 1 and 16 are excluded, and "nested parentheses" plus left-right (top-bottom) ordering give us the Sedenions listed in the first line as column and row headings, associated with the Pathions of index > 16 listed below them:

| 2 | 4 | 6 | 8 | 10 | 12 | 14 | 15 | 13 | 11 | 9 | 7 | 5 | 3 |
|---|---|---|---|----|----|----|----|----|----|---|---|---|---|
| 19| 21| 23| 25| 27 | 29 | 31 | 30 | 28 | 26 | 24| 22| 20| 18|

|   | 2 | 4 | 6 | 8 | A | C | E | F | D | B | 9 | 7 | 5 | 3 |
|---|---|---|---|---|---|---|---|---|---|---|---|---|---|---|
| **2** |   | *6* | 4 | *A* | 8 | *E* | C | *L* | F | *9* | B | *5* | 7 |   |
| **4** | *6* |   | 2 | *C* | *E* | 8 | A | *B* | *9* | F | *L* | *3* |   | 7 |
| **6** | 4 | 2 |   | E | C | A | 8 | *9* | *B* | *L* | *F* |   | *3* | *5* |
| **8** | *A* | *C* | E |   | 2 | 4 | *6* | 7 | *5* | *3* |   | *F* | *L* | B |
| **A** | 8 | *E* | C | 2 |   | *6* | 4 | *5* | 7 |   | *3* | *L* | F | *9* |
| **C** | *E* | 8 | A | 4 | *6* |   | 2 | *3* |   | 7 | *5* | *B* | *9* | F |
| **E** | C | A | 8 | *6* | 4 | 2 |   |   | *3* | *5* | 7 | *9* | *B* | *L* |
| **F** | *L* | *B* | *9* | 7 | *5* | *3* |   |   | 2 | 4 | *6* | 8 | A | C |
| **D** | F | *9* | *B* | *5* | 7 |   | *3* | 2 |   | *6* | 4 | A | 8 | *E* |
| **B** | *9* | F | *L* | *3* |   | 7 | *5* | 4 | *6* |   | 2 | C | *E* | 8 |
| **9** | B | *L* | *F* |   | *3* | *5* | 7 | *6* | 4 | 2 |   | E | *C* | *A* |
| **7** | *5* | *3* |   | *F* | *L* | *B* | *9* | 8 | A | C | E |   | 2 | 4 |
| **5** | 7 |   | *3* | *L* | F | *9* | *B* | A | 8 | *E* | *C* | 2 |   | *6* |
| **3** |   | 7 | *5* | B | *9* | F | *L* | C | *E* | 8 | *A* | 4 | *6* |   |

For compactness of presentation, edge signs associated with emanating the Assessors in the cell entries are indicated by putting their Sedenion indices in *underscored italics*, with hex notation used for indices > 9: Sedenions 10 through 15 display as letters A through F. Our table will then look as shown at left.

Cell-entry values are deployed with mirror symmetry about the void diagonals; they also show reverse-order symmetry in opposite rows: e.g., entries for the leftmost $2^{nd} - 7^{th}$ cells of the top row echo those for the $14^{th} - 8^{th}$ on the bottom-most row.

Note that all 168 available cell-entries are filled with emanations: save for the diagonals, all row-listed ZD's form DMZ's with all column-listed ZD's. As a harbinger of things to come, however, the eighth such table (strut constant indexing the Sedenions' generator) shows a peculiar alignment: in each quadrant, all edge-signs are identical. In upper left and lower right, they're negative, and positive otherwise. As is readily shown by writing products in CDP form and recalling that all inner *xor*'s will equal 16 + 8, this is a side-effect of the peculiar characteristic of generators mentioned earlier: namely, all their products with lower-indexed units, if taken on the same side, have the same sign.

More subtle, and far more surprising, is an effect ultimately attributable to the same root cause, but by a less transparent route: when the strut constant gets bigger than the generator's index, ***only 72*** of the available cell-entries contain emanations. Even more mysterious, if we subtract 8 from the strut constant and call this the "excess," then 48 of these cell entries saturate the rows and columns at a distance of the excess from the outer boundary: hence, for strut constant 9, all the outermost cells save for the 4 corners are filled, and for 15, a cross of central rows and columns (minus the central square the diagonals make) is formed, with contents evenly divided among allowed row- and column- heading Sedenions, save for '8' and the excess. Instead, these excepted entries are equally dispersed among the remaining 24 filled cells, with the edge-signs always negative for entries of index '8' (and always positive for index equal to the excess), and restricted to the same quadrants as the '/' and '\' empty diagonals respectively.

Viewed in sequence, these tables suggest the patterns made by cellular automata; seen individually, they suggest nothing so much as the short-shelf-life "sand mandalas" of Tibetan Buddhist ritual, made by monks on large flat surfaces with colored sands or powdered flowers, minerals, or even gemstones. To best convey this sensibility, all 7 of the 72-entry tables (prefaced by the "edge-polarized" table generated by strut constant 8) are arrayed in two rows of four each, to be scanned left to right, top row first. They have also been color-coded: the void diagonals are left empty; the cells which "should" be filled but are magically left content-free are painted gray; the 24 entries made of strut-constant components are green; the 48 which contain all other allowed entries are rendered in red. (Rendered originally in MS Word tables, the cell-entry contents are quite readable at 200% zoom in their native application; they become muddy, however, and hard to read, when zoomed in the PDF format they were ported to, for which we offer our apologies. In one sense, though, this "bug" is a "feature": it demands you attend to the overall patterning, which is obviously what matters most here!)

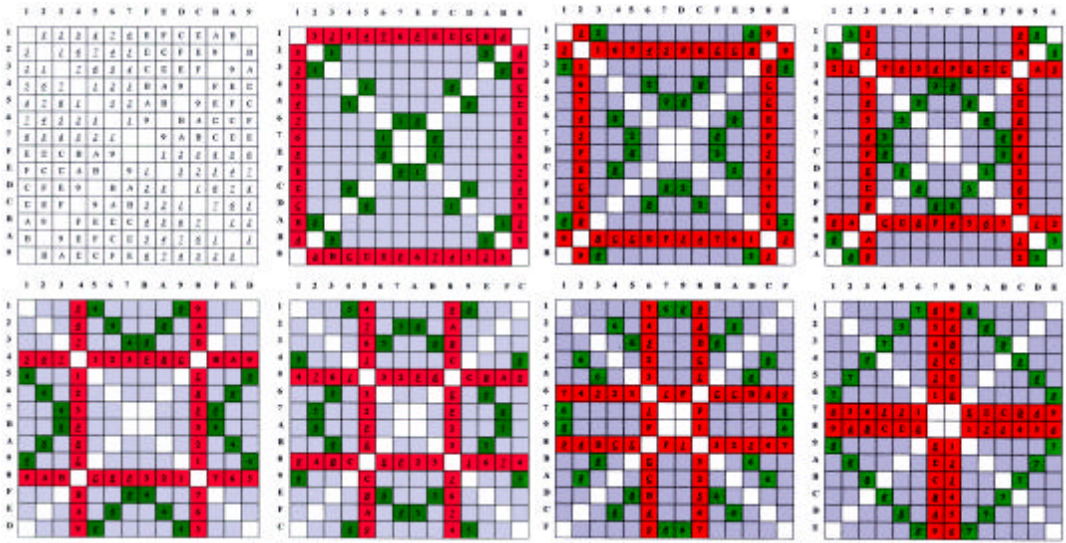

We're not yet done with surprises: for just as Kite-Chain Middens indicate a way in which low-level structure can be echoed *higher up* the endless $2^n$-ion ladder, the seven 72-cell sand mandalas show us the way to collapse from one level of form to that native to a level *lower down*. If we fold any of their tables along their row and column midlines, the red cells of the sand mandala in question will overlay precisely and only those **O**'s and **S**'s which form the vertices of the Box-Kite whose strut constant is equal to the "excess." The cells painted green, on the other hand, will adjoin the Sedenions' index-8 generators and the Box Kite's own strut constant in a trio of ensembles, one per strut – or, said another way, one for each "forbidden" **O**-triplet containing the strut constant.

We can imagine the Pathion units associated uniquely with each emanation-table entry being "boiled off" by some reductive process, the possible nature of which must be left for future contemplation. We can also, of course, imagine "Sand-Mandala-Chain Middens," echoing their basic structure in a modulo-32 rhythm. And we can imagine even richer possibilities with the 31 emanation tables associated with the Chingons . . . "restless triplet" cyclings connecting all such systems . . . and so on.

To reach an understanding of the happenings in higher reaches, a more formal grasp of the simplicities underlying the surprises in the Pathions would stand us well. We know that units indexed $g = 2^n$ – hence, any *generators* of the next ($2^n$-1) units via CDP – operate in a uniquely uniform way. For all units $i_u$, $u < g$, $i_u \cdot i_g = + i_{(u+g)}$. In the upper left ("NW") quadrant of the emanation table, meanwhile, all row and column index-pairs can be written **(r, R)** and **(c, C)** respectively, where lower-case indices all are **O**'s, and upper-case are **P**'s (all > 16). Specifically, for all 72-celled sand mandalas, each Assessor's *inner xor* ("ix") = **24** + **X**, where **X** ("xi") is the "excess" in the range of 1 to 7.

But 24 is the sum of the first two generators which yield ZD behavior. Multiplied with **O**'s, the combined effect of the unilateral signing behavior of each taken separately is *reversed*: $i_{24} \cdot i_u = + i_{(24+u)}$. Further, if $i_u \cdot i_v = + i_w$ makes an **O** triplet, $i_{(24+u)} \cdot i_{(24+v)} = + i_w$ also. (By CDP, we get $(- i_{(8+u)} \cdot i_{(8+v)})$; this gives $+ i_w$ by the "Quadrant III" rule applied to generation of **S**'s from **O**'s.) Hence, for all **u** = **r** or **c** in the NW box, the term-by-term signing of the two pairs of units in the product **(r, 24 + r\*) · (c, 24 + c\*)** will be the same as for **(r, r\*) · (c, c\*)**. But this latter just involves **O**'s, and so lacks the sign-canceling among like-indexed terms DMZ's require – and so, therefore, must the former!

The only way around this null result is for one of the indices of **R**, **C**, or **R** *xor* **C** to be *precisely* 24, with corresponding index **r**, **c**, or **r** *xor* **c** precisely the excess **X**. Say **r** = **X**, left-multiply **(c, C)** by **(r, R)**, and leave a sign-variable **sg(..)** = ± **(..)** undetermined:

**(X, 24) · (c, 24 + c\*)**, where **c · c\*** = sg( **X** ), and hence **X · c** = sg( **c\*** )

Left-multiplying by **X**, then **24**, we get these indices in successive output lines:

**X · c** → sg( **c\*** ),     **X · (24 + c\*)** → 24 + (**X · c\***) = ( – sg(24 + c) )
**24 · c** → (24 + c),     **24 · (24 + c\*)** → ( – **c\*** )

Comparing cross-terms in the two lines, it's clear we have a zero product, with the edge between **(r, R)** and **(c, C)** being signed the same as **sg**. Supposing **c** is **X** instead leads to the same result with a change of labels; slightly more work is called for if **X** is in the *product* instead of either multiplier. By symmetry, the same ZD criteria just demonstrated can be appropriately interpreted in all quadrants. The 72-count (not just of non-void cell entries, but of distinct DMZ pairings of ZD diagonals) follows trivially.

That last parenthetical remark is a signal one: in Sedenion Box-Kites, the 12 ZD diagonals in each, each belong to 2 Sails, making 24 DMZ's in all. In this, they conform to the same formula we can generalize from the sand mandalas with **X** > 0: for $2^n$-ions, n>3, the ZD-pairings formula, up to Voudons, gives $6 \cdot (2^{n-1}-4)$ = 24, 72, 168, 360, 744. Alone in this series, the pairings for Pathions are a simple multiple of their predecessors'. A little investigation shows us why: each of the 7 "excessive" sand mandalas can in fact be partitioned into 3 Box-Kites, all sharing the same strut constant, but otherwise exactly like those in the Sedenions – save for their each containing an "impossible" 8 as the low index in one Assessor, forming a triple zigzag with two Assessors whose low-index terms are an **O** and its *xor* with 8! So not only can we "fold over" its emanation table to get a Sedenion Box-Kite; we can see it as housing 3 "harmonics" of it of a special kind. The analogy to the Hyperbolic vs. Elliptic Umbilic Catastrophes' 1 to 3 Cusps, *vis à vis* 1 triple-zigzag to 3 trefoil Sails, was explored in our prior paper. How deep does it go?

In Catastrophe Theory, the $D_4$ singularity takes on two forms when projected from the complex to the real domains, and these are the two Umbilics just named. More, it is also well-known, thanks to the "A-D-E Problem,"[20] that the $D_4$ pattern relates to close-packing in 4-D, yielding the "Feynman checkerboard" of 24 unit hyperspheres in each other's immediate proximity (the same pattern giving the 24-element group of integer Quaternions). The Hyperbolic Umbilic, meanwhile, has its 2-D "behavior plane" – when its definitive control variable is zero – degenerate into the 4-to-1 collapsing of quadrants topologists call "folding the handkerchief." The similarity of this process to the "folding over" of emanation tables – plus the special significance of 24 here – seems hardly accidental. (And 72, for what it's worth, is the "closest-packing number" associated with the singularity $E_6$ – the basis of the mysterious "trance tunnel" geometry studied by James Callahan in a series of seminal papers modeling emotional affect response.[21]) How might these patterns echo in higher number forms? We have not, for instance, considered the "pathologies" awaiting in higher $2^n$-ions, when not 2, but 3 or more such "hypernumbers" are multiplied. (Will we need, perhaps, to contemplate "emanation *solids*"? And does the fact that both the Sedenions' 24 and the Pathions' 72 are contained as simple factors of the Routons' "Babylonian" 360 mean anything special? In this domain, so little is known, and so much computation will be required to learn more, that we can truly deem this brand of mathematics an "experimental" discipline.)

Let's wind down with something more concrete to mull on: the emanation table shown as the last on the right of the top row in the prior illustration (with **X** = 3) will be exemplified. Its "green cell" units of 3 and 8 indicate a "folding over" into the Sedenion Box-Kite with strut constant 3; the low-end indices shown in the emanation table do indeed "fold over" in the manner described above so that strut-opposed terms form Sedenion Assessors (2, 9), (1, 10), (7, 12), (6, 13), (5,14), (4,15) as required. Meanwhile, if you check the first table presented in this paper to corroborate that statement, you'll recall the format of presentation (Assessors A, B, C forming a triple-zigzag, etc.). That format can be retrieved (*sans* the now-redundant strut constant in the left-most column) to present the 3 "impossible" Box-Kites the emanation table in question partitions into (each of whose triple-zigzags' bolded lower-index terms form hard-working, Box-Kite-sustaining "mule-trains" of "sterile" triplets which can never form ZD's in their own right). This emanation table partitions into 3 Box-Kites with strut constant 11 as follows:

| A | B | C | D | E | F |
| --- | --- | --- | --- | --- | --- |
| **1**, 26 | **8**, 19 | **9**, 18 | 2, 25 | 3, 24 | 10, 17 |
| **4**, 31 | **8**, 19 | **12**, 23 | 7, 28 | 3, 24 | 15, 20 |
| **6**, 29 | **8**, 19 | **14**, 21 | 5, 30 | 3, 24 | 13, 22 |

Note that the common appearance of the (8, 19) and (3, 24) Assessors in each of these 3 Box-Kites would seem to belie the assertion that they partition the 72 DMZ cells cleanly: but this shared pair of Assessors are strut-opposites, and we are now viewing these Box-Kites *not* as collections of 12 ZD diagonals, but of 24 *DMZ-pairings*; so, *as such*, there's *no sharing* of any terms between them! From all we've seen of the Pathions and their patterns of emanation, this distinction, while subtle, is obviously a crucial one.

But this suggests the 7 interconnected Box-Kites of the 168-cell emanation tables may be viewed in a similar manner, making their interconnectedness not hard to deal with, but simple. Consider the first table's Assessor list given atop page 14: take the lower-index items in the first of its two lines as given, and pick the leftmost value (2). The two successive pairs that follow each complete a triple with it. Add in their 3 strut-opposites on the right, and we've got the lower-index scheme for the first 3 Box-Kites: (**2, 4, 6**, 7, 5, 3); (**2, 8, 10**, 11, 9, 3); (**2, 12, 14**, 15, 13, 3). Now move right, pick the (4), and form 2 new triples with it from the left of the line; including opposites from the right gives 2 more Box-Kites: (**4, 8, 12**, 13, 9, 5); (**4, 10, 14**, 15, 11, 5). Shift again, pick the (6), perform the same trick: (**6, 8, 14**, 15, 9, 7); (**6, 10, 12**, 13, 11, 7).

Parting thought anent "Mule-Trains": the attentive reader may already have figured out that if each **X**-high (henceforth, "Sky-High") emanation table splits into 3 Box-Kites, each with 2 sails based on one of the 7 overworked "Mule-Train" triplets, 7 x 3 x 2 = 42. The MT-Sails, then, play a role in Sky-Highs akin to the Sedenions' 42 Assessors – or, perhaps more precisely, to the 42 inscriptions of the similarly "sterile" Octonions on Box-Kite vertices therein. In this sense, each Sky-High's excluding *but one* Mule-Train (the triplet completing the "green cell" pairings) is like the ban of the strut-constant *Octonions* in the Sedenion context, which also act as "signatures." Is all this the beginning of a "Jordan-Holder" sort of theory of lattice hierarchies in the ZD-Net? Tune in next monograph, same text archive, same font of speculation, for the next exciting episode . . .

**{3} Simulating "Creation Pressure" ( Cellular Automata Underneath Number?):** In the early days of neural networks, random Booleans operator values were strewn across gate arrays, while researchers "watched what happened." Differently weighting different levels of $2^n$-ions suggests itself, with further study of the harmonics and collapses, first in evidence with the Pathions as just discussed, calling for a more programmatic approach. The possibility of underwriting such a "ZD-net" with cellular-automata-like rules – and hence, seeing Number itself as their side-effect in some deep senses – merits serious consideration as well. Steve Wolfram's "Rule 30,"[22] for instance – if a term and its right-hand neighbor share a property, they switch its value to that held by the neighbor on their left – might prove a good starting point in exploring "triplet-switch effects": for any Assessor, in any Sail, its lower-order (**O**, say) index will form a triplet with the other lower-order indices in the same Sail; but it will also form a triplet with its adjacent Assessors' higher-order (**S** or **P**, say) indices.

Meanwhile, collapsing a Sky-High onto a Sedenion box-kite leaves all Assessors properly indexed, but since the strut-opposite cell contents from which they were engendered in the Pathions always have opposite edge-signs, their *own* edge-signs are left indefinite, and hence any interrupted cyclings may perhaps be left hanging (and vulnerable, thereby, to "digression pressure") . . . at least for the moment of symmetry-breaking when Pathions are "boiled off."

Then too, kite-chain midden harmonics allow for ever greater "switching yard" complexities as the n in $2^n$-ion gets bigger. What underlying restlessness of triplets, or overarching "gantry-work" for growing n, might lead all these systems to "cross-pollinate"? To what effect? By induction, and in conclusion, one can easily imagine things in the general case getting wacky with alacrity. Where, then, one must wonder, might we go with this?

[1] R. Guillermo Moreno, "The zero divisors of the Cayley-Dickson algebras over the real numbers," Bol. Soc. Mat. Mexicana (3), 4, 1 (1998), 13-28; preprint available as http://arXiv.org/abs/q-alg/9710013 .

[2] K. and Mari Imaeda, "Sedenions: Algebra and analysis," Appl. Math. Comput., 115, 2/3 (October 2000), 77-88.

[3] Robert de Marrais, "The 42 Assessors and the Box-Kites they fly: Diagonal axis-pair systems of zero-divisors in the Sedenions' 16 dimensions," preprint available as http://arXiv.org/abs/math.GM/0011260 . (November, 2000)

[4] Corrado Segrè, in Math. Ann., 40 (1892), 413-167, and in Segrè's Opere, I, Edizione cremonese, 1957.

[5] G. Baley Price, Introduction to Multicomplex Spaces and Functions, Marcel Dekker, New York, 1991.

[6] R. Penrose and W. Rindler, Spinors and Space-Time, Vol. 1: Two-Spinor Calculus and Relativistic Fields, Cambridge University Press, Cambridge – London – New York, 1984.

[7] Charles Musès, "Applied Hypernumbers: Computational Concepts," Appl. Math. Comput., 3 (1976), 211-226; "Hypernumbers – II. Further Concepts and Computational Applications," Appl. Math. Comput., 4 (1978), 45-66; "Hypernumbers and Quantum Field Theory with a Summary of Physically Applicable Hypernumber Arithmetics and their Geometries," Appl. Math. Comput., 6 (1980), 63-94.

[8] Pertti Lounesto, Clifford Algebras and Spinors, Cambridge University Press, Cambridge – New York – Melbourne, 1997, 207-8.

[9] Andrei Sakharov, "Vacuum quantum fluctuations in curved space and the theory of gravitation," Dokl. Acad. Nauk. SSSR 177 (1967), 70-71; translated in Sov. Phys. Dokl. 12 (1968), 1040-1041.

[10] H. E. Putoff, Phys. Review A, 39 (1989), 2223; B. Haisch, A. Rueda, H. E. Putoff, Phys. Review A, 49 (1994), 678; B. Haisch, A. Rueda, Phys. Review A, 268 (2000), 224; Y. Dobyns, A. Rueda, B. Haisch, xxx.lanl.gov/gr_qc/0002069 (February, 2000), to appear in Foundations of Physics.

[11] Jaak Lõhmus, Eugene Paal, Leo Sorgsepp, Nonassociative Algebras in Physics, Hadronic Press, Palm Harbor FL, 1994. Geoffrey Dixon's extension of the Standard Model based on the tensor product of the Complex, Quaternion and Octonion algebras, where electric, weak, and strong forces naturally arise from these 3 factors respectively, suggests the next "big break" – allowing zero-divisors – would point toward including gravity (and that nothing *less* than this could hope to). See Dixon's Division Algebras: Octonions, Quaternions, Complex Numbers and the Algebraic Design of Physics, Kluwer, Dordrecht, 1994.

[12] John Milnor, Morse Theory, based on lecture notes by M. Spivak and R. Wells, Princeton University Press, Princeton NJ, 1969. Diagram and quote from p. 54, at the end of "§9. The Curvature Tensor."

[13] J. J. Sylvester, "On quaternions, nonions, sedenions etc.", Johns Hopkins Univer. Circular, 3:7-9, 1889; this and earlier pieces on Nonions are in his Collected Works, III, 647-650; and IV, 122-132.

[14] D. Siersma, "Singularities of $C^8$ functions of right-codimension smaller or equal than eight," Indag. Math. 74 (1973), 31-37; simply explained and motivated in Tim Poston and Ian Stewart, Catastrophe Theory and its Applications, Pitman Publishing, Boston – London – Melbourne, 1978, pp. 162-170.

[15] A. R. Rajwade, Squares, Cambridge University Press, Cambridge – New York – Melbourne, 1993.

[16] T. Kirkman, "On pluquaternions and homoid products of sums of *n* squares," Philos. Mag. (Ser. 3) 33 (1848), 447-459; 494-509.

[17] I've yet to find the actual paper, "The algebraic set of the zero divisors in the Cayley-Dickson algebras," but its results are cited at http://www.innerx.net/~tsmith/Ndalg.html -- Tony Smith's well-traveled math and physics website. A brief abstract of Moreno's talk, in which he presents the results cited above as "*upper bounds*" on the number of "irreducible components" (a point of much interest to us) can be found at www.matmor.unam.mx/events/ams-smm/EN/SESIONES/sessionR1/sesionR1.pdf .

[18] Helena Albuquerque and Shahn Majid, "Quasialgebra structure of the octonions," preprint available as http://arXiv.org/abs/math.QA/9802116 . (February, 1998)

[19] Moreno's formula is identical to that derived by Hurwitz in a superficially unrelated context: the classification of compact Riemann surfaces of genus $g \geq 2$ have automorphism groups of size 84(g-1). And, as the smallest order such group is a simple group, there is no Hurwitz group of order 84 – so the first instance is analogous to the Pathion harmonic just described, while higher instances are no longer general! (But indeed, even the Pathions' case isn't truly general either, as we'll see next section!) See Gareth A. Jones and David Singerman, Complex Functions: An Algebraic and Geometric Viewpoint, Cambridge University Press, Cambridge – New York – Melbourne, 1987, Chapters 5 and 6.

[20] M. Hazewinkel, W. Hesselink, D. Siersma, F. D. Veldkamp, "The ubiquity of Coxeter-Dynkin diagrams (an introduction to the A-D-E problem), Nieuw Arch. Wisk. 25 (1977), 255-307. The A-D-E problem arose in the work of V. I. Arnol'd as a result of his sensing the connection between singularity theory and

the classic "A, D, E" of Dynkin diagrams. See his Catastrophe Theory (3$^{nd}$ revised and expanded edition),Springer-Verlag, Berlin Heidelberg New York Toronto, 1992, especially the concluding chapter. An updated roster of links on the subject is maintained by Tony Smith at the website cited earlier.

[21] James Callahan, "A Geometric Model of Anorexia and Its Treatment," Behavioral Science, 27, 1982, 140-154, is where the $E_6$ dynamics of the "trance tunnel" were first reported; these were embedded in the full dimensionality of the Double Cusp in a pair of papers written with psychoanalyst Jerome I. Sashin, "Models of Affect-Response and Anorexia Nervosa," in S. H. Koslow, A. J. Mandell, M. F. Shlesinger, eds., Conference on Perspectives in Biological Dynamics and Theoretical Medicine. Annals N. Y. Acad. Sci., 504, 1986, 241-259; "Predictive Models in Psychoanalysis," Behavioral Science, 35, 1990, 60-76.

[22] Steve Wolfram, A New Kind of Science, Wolfram Media Inc., www.stevewolfram.com, 2002.